\documentclass[]{article}

\title{Computer Validation of Neural Network Dynamics:\\ A First Case Study}
\author{Christian Kuehn\thanks{Department of Mathematics, Technical University of Munich, 85748 Garching b. München, Germany}~\thanks{Munich Data Science Institute (MDSI), Technical University of Munich, 85748 Garching b. München, Germany}~~and Elena Queirolo\thanks{Department of Mathematics, Technical University of Munich, 85748 Garching b. München, Germany}}
\date{}
\usepackage{mystyle}

\begin{document}

\maketitle

\begin{abstract}
A large number of current machine learning methods rely upon deep neural networks. Yet, viewing neural networks as nonlinear dynamical systems, it becomes quickly apparent that mathematically rigorously establishing certain patterns generated by the nodes in the network is extremely difficult. Indeed, it is well-understood in the nonlinear dynamics of complex systems that, even in low-dimensional models, analytical techniques rooted in pencil-and-paper approaches frequently reach their limits. In this work, we propose a completely different perspective via the paradigm of validated numerical methods of nonlinear dynamics. The idea is to use computer-assisted proofs to validate mathematically the existence of nonlinear patterns in neural networks. As a case study, we consider a class of recurrent neural networks, where we prove via computer assistance the existence of several hundred Hopf bifurcation points, their non-degeneracy, and hence also the existence of several hundred periodic orbits. Our paradigm has the capability to rigorously verify complex nonlinear behaviour of neural networks, which provides a first step to explain the full abilities, as well as potential sensitivities, of machine learning methods via computer-assisted proofs. We showcase how validated numerical techniques can shed light on the internal working of recurrent neural networks (RNNs). For this, proofs of Hopf bifurcations are a first step towards an integration of dynamical system theory in practical application of RNNs, by proving the existence of periodic orbits in a variety of settings.
\end{abstract}


\section{Introduction}

Machine learning~\cite{Alpaydin,Murphy,SchoelkopfSmola} has been an extremely dynamic field in recent years~\cite{JordanMitchell,Silveretal}. Deep neural networks have taken center stage as a key algorithmic component~\cite{LeCunBengioHinton,Schmidhuber}. Thanks to the results presented in \cite{chen2018neural}, and subsequent developments, it is now possible to explicitly connect several neural network architectures with an explicit dynamical system given by ordinarny differential equations (ODEs).
Although several results about static functional approximation properties of neural networks are quite classical~\cite{Pinkus,ScarselloTsoi}, studying their dynamics~mathematically rigorously, i.e., to really \emph{prove} that they generate a particular behaviour, still remains a formidable challenge. From an abstract viewpoint, this is not too surprising since the propagation of information across a neural network with fixed edge weights, as well as the learning process of the edge weights, are large-scale nonlinear dynamical systems~\cite{E5,FarmerPackardPerelson}. It is well-understood that if one seeks rigorous mathematical results in nonlinear dynamics, one usually has to make quite significant assumptions~\cite{GH}. Examples of these assumptions are low dimensionality~\cite{DeMeloVanStrien,Sanger}, scale separation~\cite{GaltierWainrib,KuehnBook}, existence of a mean-field~\cite{MeiMontanariNguyen,Tanaka}, or a gradient structure~\cite{AmbrosioGigliSavare,WangTengPerdikaris}. In fact, if these assumptions are missing, then already questions about very low-dimensional dynamics, say in dimensions one to three for iterated maps or flows, become almost impossible to solve via rigorous pencil-and-paper arguments~\cite{Ilyashenko2}. Therefore, if we are interested in generic neural networks without simplifying assumptions, which are frequently the cases most relevant for many practical algorithms, it seems that we would have to abandon the hope to ever validate mathematically that modern machine learning algorithms behave correctly. 

In this work, we propose a completely new perspective to validation of neural network dynamics and thereby also to robustness and verification of artificial intelligence. The idea is to use automated computational tools to rigorously validate the dynamics of neural networks. Although this perspective might sound counter-intuitive at first, it can be well-founded within the field of numerical validation methods~\cite{ArioliKoch1,DayLessardMischaikow,Kapelaetal,vandenBergJamesReinhardt}. In particular, viewing the neural network as a dynamical system, we are going to utilize and adapt computational techniques first developed to solve hard open problems in dynamics, probably most famously the existence of the Lorenz attractor~\cite{Tucker}. The principle is that certain patterns, respectively solution behaviours, can be tackled numerically. The numerical solution is then verified rigorously via interval arithmetic error bounds in combination with a-priori or a-posteriori error estimates. This combination of techniques can indeed provide validated mathematical proofs for otherwise intractable global nonlinear dynamics. In the context of neural networks, it has the nice interpretation that we use a more classical computational approach in combination with theoretical mathematical bounds and computer-assisted error tracking to validate the existence of patterns in neural networks, which can then be used as an even more powerful computational tool.  

Of course, implementing such an idea requires a longer term development and testing for various classes of deep neural networks. In this work, we are interested in a first proof-of-principle, i.e., to demonstrate how the idea of numerical validation can be brought to tackle intermediate size neural networks. Both the method and the code presented here are given in great generality and 
as a test case, we study a well-established class of anti-symmetric recurrent neural networks (RNNs)~\cite{chang2019antisymmetricrnn}. In this setup, we focus on the dynamics on the network for edge weight matrices with random weights to ensure robustness of our results. We show that in a generic setting and under the assumption of large depth, varying a main hyperparameter of the network leads to a whole finite cascade of Hopf bifurcations, where periodic solutions are generated. Via the so-called \emph{radii polynomial  approach}, we validate several hundred Hopf bifurcations rigorously and we also compute the relevant first Lyapunov coefficient rigorously. We find that the generated periodic solutions are all unstable leading to a phase space structure for the neural network, where many transient oscillatory motions are possible. In fact, it is well-understood that within the class of all possible dynamical systems, looking for large classes of (unstable) periodic solutions is a strategy to determine the possible high complexity and potentially chaotic dynamics of the system. In summary, we have demonstrated rigorously that within the space of even varying a single hyperparameter, the considered class of RNNs has enormous dynamical complexity. 
Similar complexity is retrieved in a variety of other networks, as shown in in Section \ref{sec:examples}.
More importantly, we have demonstrated a paradigm that can complement validating the dynamics of neural networks precisely in the regimes, where other, often quite strong, mathematical assumptions needed for purely analytical arguments break down.

The paper is structured as follows: in Section \ref{sec:finding_Hopf} an initial presentation of Hopf bifurcations is given, building towards the algebraic problem defining them. Then, Section \ref{sec:RPA} presents a bird eye view of validated =numerics in the context of finite dimensional problems, such as the Hopf algebraic equation. These two ingredients are the mathematical background to then discuss in depth our first example: Antisymmetric RNNs, presented in Section \ref{sec:BfInRNN}. There we give first a justification on the choice of network and of bifurcation, then present extensive results on the effectiveness of our validation method. More general examples are provided in Section \ref{sec:examples}, where a wide variety of network structures and continuation parameters are showcased. Additional examples can be found directly in the code at \cite{ODEinRNN}. Finally, Section \ref{sec:outlook} provides a brief  outlook. Two Appendices are given with details on the application of the radii polynomial approach to the RNN (Appendix \ref{sec:RPA_details}) and on the computation of the Lyapunov coefficient to prove non-degeneracy (Appendix \ref{sec:FLC}).

\section{Algebraic Hopf bifurcation}
\label{sec:finding_Hopf}

Following \cite{chen2018neural}, RNNs can be either interpreted or directly coded as a system of ordinary differential equations called \emph{neuralODE}. As an overview, when considering RNNs their formulation is 
\begin{equation}\label{eq:discreteRNN}
    x_{t+1} = x_{t} + \sigma(x_t, \theta),\quad t=1,2,\dots, N
\end{equation}
where $x_0$ is the input layer, $x_N$ is the output layer, $t$ is the hidden layer, $\sigma$ is the chosen non-linearity and $\theta$ are all the parameters optimized during the learning process, usually it included the weights and biases. Assuming that $\sigma$ can be rewritten as $\epsilon f$ for a small $\epsilon$, then we can recognise in Equation \eqref{eq:discreteRNN} an application of the forward Euler method to the ODE
$$ 
x'(t) = f(x(t),\theta), \quad t\in[0, T],
$$
where $T=\epsilon N$. From now on, $t$ is considered to be the continuous variable unless otherwise specified. We are going to use approaches and techniques rooted in dynamical systems to gain an understanding of the behaviour of RNNs in the limit of neural ODEs; yet, it is important to point out that the principles we employ will also be useful for neural networks beyond ODEs. In this paper we will concentrate on Hopf bifurcations, considered as places of onset of periodic orbits. To further justify this choice in a particular setting, in the following Section \ref{sec:BfInRNN}, the possibility of having Hopf bifurcations appearing in the structured RNN \eqref{eq:ODE_AsymRNN} is analytically justified. 

We now discuss the numerical search for such bifurcations, by introducing the zero finding problem associated to Hopf bifurcations. 
Since Hopf bifurcations depend on only one parameter, we consider here the situation in which $\theta = \gamma\in\mathbb{R}$, where $\gamma $ is the unique parameter of the system. In the future, when considering full RNNs with more than one parameter, we will assume only one parameter is flexible, while all others are fixed.

\begin{definition}
\label{def:Hopf}
A Hopf bifurcation $(x_\star, \gamma_\star)$  of $x' = f(x, \gamma)$ is such that $(x_\star, \gamma_\star)$ is an equilibrium  of the ODE, i.e. $f (x_\star, \gamma_\star)= 0$, such that $\textnormal{D}_x f (x_\star, \gamma_\star)$ has a pair of purely imaginary eigenvalues $\lambda_\star, \bar\lambda_\star$. Furthermore, the Hopf bifurcation is non-degenerate if there are no other imaginary eigenvalues of $\textnormal{D}_x f (x_\star, \gamma_\star)$ and $\lambda_\star$ crosses the imaginary axis with non-zero velocity w.r.t.~$\gamma$ at $\gamma = \gamma_\star$.
\end{definition}

Setting momentarily aside the non-degeneracy conditions, we can set the algebraic problem as
\begin{align*}
\begin{cases}
&f(x, \gamma) = 0,\\
&\textnormal{D}_xf(x, \gamma) v - \lambda v = 0, \qquad \lambda \in \imag \mathbb{R}.
\end{cases}
\end{align*}
Notice how the eigenvalue $v$ is defined just up to a scaling, thus an additional equation involving the normalization of $v$ needs to be included to guarantee uniqueness of the solution. We define all solutions of this problem as \textit{algebraic Hopf bifurcations}. Such bifurcations might be degenerate.

While the  algebraic Hopf problem is mathematically well-posed, up to a complex rescaling of $v$, it is defined on the space $(x,\gamma, v, \lambda) \in (\mathbb{R}^n, \mathbb{R},\mathbb{C}^n, \imag\mathbb{R} )$. It is numerically cumbersome to impose that the solution exactly fits into the appropriate space, since it mixes real, complex and imaginary values. Instead of solving it directly, we rephrase the problem  into a fully real space. We introduce the notation $v = v_r + \imag v_i$ and $\lambda = \imag \lambda_i$, where now $v_r, v_i\in \mathbb{R}^n$ and $\lambda_i\in\mathbb{R}$. 
We also add two equations
$$
\phi^\top v_r = 0 \qquad \phi^\top v_i  - 1 = 0,
$$
where $\phi$ is any fixed vector in $\mathbb{R}^n$, to fix the scaling of $v_r $ and $ v_i$. Then, we write the real algebraic Hopf problem as
\begin{equation}\label{eq:Hopf}
\begin{cases}
&\phi^\top v_r = 0, \\ 
&\phi^\top v_i  - 1 = 0,\\
&f(x, \gamma) = 0,\\
&D_x f(x, \gamma) v_r + \lambda_i v_i =0 , \\
&D_x f(x, \gamma) v_i - \lambda_i v_r = 0 . \\
\end{cases}
\end{equation}
With the addition of this scaling, the problem is numerically well-posed and has, generically,  a locally unique solution. Equation \eqref{eq:Hopf} is now real and finite dimensional, any root finding algorithm, such as Newton's method, can be used to find its numerical solutions. This yields numerical approximations of algebraic Hopf bifurcations. A discussion on the initialization of such root finding algorithms can be found in Remark \ref{rem:initial}.

Having computed a numerical solution to the algebraic Hopf problem, an overview of its validation is given in Section \ref{sec:RPA}, while details pertaining to its practical implementation are presented in Appendix \ref{sec:RPA_details}. The problem of proving non-degeneracy is discussed at the end of the Section \ref{sec:RPA}, while details of its computation are presented in Appendix \ref{sec:FLC}.

\section{Validation in finite dimensions}\label{sec:RPA}

In this section, an \emph{a-posteriori} method of validation for zero-finding problems is presented. Following, among others, \cite{hungria2016rigorous, lessard2016automatic}, we will give here an overview of the radii polynomial approach. 

Consider \eqref{eq:Hopf} as a zero-finding problem $F(x)=0$, where $F:X\rightarrow Y$ and $X,Y$ are Banach spaces with norms $\|\cdot\|_X$, $\|\cdot\|_Y$. For concreteness, one may think of finite-dimensional Euclidean spaces $X,Y$ here but the validation idea works in more generality, even in infinite-dimensional settings, so we keep this generality in the presentation to make it evident, how far-reaching the approach actually is. If $\tilde x$ is close enough to a solution of the zero-finding problem, we expect the Newton operator $\tilde x - \textnormal{D}F(\tilde x)^{-1} F(\tilde x)$ to be contracting towards the exact solution. Based on this intuition, we define  the map
\begin{align}\label{eq:contraction_mapping}
T: X&\rightarrow X \nonumber\\
 x &\mapsto x - A F(x),
\end{align}
where $A: Y \rightarrow X$ is an approximation of $DF(\tilde x)^{-1}$. We the set out to prove that $T$ is a contraction. More precisely,
with the radii polynomial approach, we prove the existence of an $r$ such that the approximate Newton operator $T$ is a contraction in the ball $B_r(\tilde x)$ of radius $r$ around the numerical solution $\tilde x$. This yields a rigorous existence of a zero very close to the numerically computed one $\tilde x$. With this strategy in place, we can now present the radii polynomial theorem from \cite{breden2013global}.

\begin{theorem}
Let $T$ be as defined in \eqref{eq:contraction_mapping} and let
$$
Y \geq \| T(\tilde x)\|_X, $$ $$ Z(r) \geq \sup_{b,c\in B_1(0)\subset X}\| \textnormal{D}T(\tilde x+rb)rc \|_X.
$$
The radii polynomial is
$$
p(r)\bydef Y + Z(r) - r. 
$$
If there is a $r_\star>0$ such that
$
p(r_\star) = Y + Z(r) - r < 0
$ 
then $T$ is a uniform contraction on $B_{r_\star}(\tilde x)$. If $A$ is non-singular, then $F$ has a unique zero in $B_r(\tilde x)$
\end{theorem}

The application to our finite dimensional system \eqref{eq:Hopf} is presented in the Appendix \ref{sec:RPA_details}. Since the bounds need to be themselves rigorously computed, interval arithmetic is necessary. 

\begin{remark}
To ensure the correctness of the bounds, all numerical errors need to be considered, such as rounding errors in the computation of the hyperbolic tangent and floating point errors. Interval arithmetic is the tool used to keep track of such errors, the used implementation in Matlab is the Intlab library,~\cite{rump1999intlab}.
\end{remark}

Once the solution to \eqref{eq:Hopf} is validated, we have proven the existence of an algebraic Hopf bifurcation, but we do not yet have knowledge of its non-degeneracy. For a non-degenerate bifurcation, we need to satisfy the two non-degeneracy conditions, as in Definition \ref{def:Hopf}. The first one is the lack of other imaginary eigenvalues. Considering a finite dimensional system, this condition can be checked directly, by computing all other eigenvalues and confirming that their real part is non-zero. Validation of eigenvalues is a built-in functionality in Intlab and to prove this condition is straightforward.

The second condition for non-degeneracy is for the imaginary eigenvalue pair to be crossing the imaginary axes with non-zero velocity w.r.t.~the parameter. This condition is more technical and it is equivalent to the first Lyapunov coefficient being non-zero, as presented in \cite{kuznetsov2013elements}. To discuss the computation and validation of the first Lyapunov coefficient in our situation, we refer the interested reader to Appendix \ref{sec:FLC}.

\section{First Example: AntisymmetricRNN}\label{sec:BfInRNN}


Recurrent neural networks are build to integrate time-dependent data in their computation, while the hidden layers keep track of previous information. 
The output layer $x_T\in \mathbb{R}^n$ is the result of applying a sequence of nonlinearities to the input layer $x_0\in\mathbb{R}^n$ accoding to
$$
x_t = x_{t-1} + f(x_{t-1}, d_t), \quad t=0,1,\dots, T
$$
where $x_t\in\mathbb{R}^n$ is the value of the $t$-th hidden layer, and $d_t\in\mathbb{R}^m$ is input for the same layer, that is additional information available at time $t$. In this formulation, $f$ is a function dependent on parameters optimized during training. The output of the network is then a time dependent value $y_t$ defined as $y_t = g(x_t)$, for some function $g$. Let us remark here that in practice neither $f$ nor $g$ need to be differentiable, but continuity is necessary. In this paper, we only consider a special case of $f$, where it is sufficiently smooth.

We are focusing on a class of neural networks known as antisymmetric recurrent neural networks (AntisymmetricRNN~\cite{chang2019antisymmetricrnn}). 
In this class of recurrent neural networks, each hidden layer is determined by
\begin{align*}
x_t = x_{t-1} + \sigma ( W x_{t-1} + V d_t + b),  \quad t=0,1,\dots, T
\end{align*}
where $W\in\mathbb{R}^{n\times n}$,  $V\in\mathbb{R}^{n\times m}$ and $b\in\mathbb{R}^n$ are parameters of the RNN. In practice, these parameters can depend on the layer, but in our exposition they will be fixed w.r.t. $t$. Furthermore, we focus on the information propagation dynamics on the trained network and not on the learning step. $W$ and $V$ are called weights, while $b$ is the bias. In specific applications, it is possible to vary the dimensions $n$ and $m$ along the layers, but we will not make use of such flexibility for this first case study. In~\cite{chang2019antisymmetricrnn}, a concrete subclass of RNNs is studied given by
\begin{equation}
\label{eq:net1}
x_t = x_{t-1} + \varepsilon \tanh (\hat W x_{t-1} + V d_t + b),
\end{equation}
where 
\begin{equation}\label{eq:Wstructure}
\hat W \bydef W- W^\top + \gamma \textnormal{Id},
\end{equation}
and
$\varepsilon\in\mathbb{R}$ and $\gamma\in\mathbb{R}$ are hyperparameters, and $\textnormal{Id}$ is the identity matrix. Here the hyperbolic tangent acts  element-wise. Following \cite{chen2018neural} one may view equation~\eqref{eq:net1} as the application of the Euler method to the ordinary differential equation (ODE) 
\begin{equation}\label{eq:ODE_AsymRNN}
x'(t) = \tanh( \hat W x(t) + V d(t) + b),
\end{equation}
where we now consider $x$ to be an $n$-dimensional function of continuous time $t\in[0, T]$. The value of $x_t$ at the  input and output layers of the RNN are equivalent to the values of $x(t)$ at $t=0$ and $t=T$. 
\begin{remark}
This rewriting is supported by a variety of results, such as \cite{zhang2019approximation,yan2019robustness}. It is possible to train a network as a neural ODE, and to \emph{a priori} chose different numerical integrators.
\end{remark}

Studying the dynamics induced by \eqref{eq:ODE_AsymRNN} is a powerful tool   in   understanding  the effects the RNN has on its input vector. Bifurcations provide insight on the global behaviour of a system through a local study, so they are a fitting starting topic. We are particularly interested in studying the bifurcations this dynamical system can exhibit depending on the hyperparameters after the weights and bias have been determined during the learning phase of the RNN. In particular, following the presentation in \cite{chang2019antisymmetricrnn}, we want to determine when the system can undergo a Hopf bifurcation. Hopf bifurcations are one of the mechanisms for the appearance of periodic orbits, a topic of interest in its own right. Hopf bifurcations happen when a pair of imaginary eigenvalues of the Jacobian crosses the imaginary axes with non-zero velocity. We refer to Definition \ref{def:Hopf} for additional details on their mathematical structure. In this section, we first justify that $\gamma$ is the hyperparameter of interest, when searching for Hopf bifurcations in  \eqref{eq:ODE_AsymRNN}.

In order to simplify the explanation, we consider no input data over time, i.e. $d(t) =0$. The unique zero of the hyperbolic tangent is at zero, so the non-linearity simplifies and we retrieve the algebraic problem
$$
\hat W x(t) + b = 0.
$$
For simplicity, and without loss of generality, we restrict ourselves additionally by setting the bias to zero, $b=0$. In this case, the unique equilibrium of the $n$-dimensional system is the trivial equilibrium $x_\star = 0 $.
The Jacobian of \eqref{eq:ODE_AsymRNN} is
\begin{equation}
\label{eq:Jacobian}
J(x) = (\textnormal{Id} - \diag(\tanh( \hat W x(t) )^2)) \hat W,
\end{equation}
where the power is considered element-wise and $\diag(\cdot)$ is the diagonal matrix obtained from a vector. Computing the Jacobian at $x_\star$, the hyperbolic tangent vanishes, leaving us with
$$
J(x_\star) = \hat W.
$$
Thus, studying the behaviour of the Jacobian at the equilibrium is equivalent to studying $\hat W$. We notice that the eigenvalues of anti-symmetric matrices, such as $W-W^\top$ are imaginary. Since $\hat W$ is a shift by $\gamma$ of an anti-symmetric matrix, the real part of all eigenvalues of $\hat W$ is equal to $\gamma$.
In particular,  all eigenvalues of $\hat W$ cross the imaginary axes for $\gamma$ crossing zero. This creates a highly degenerate Hopf bifurcation at $\gamma=0$. 

\begin{remark}
Imaginary eigenvalues appear in pairs, thus antisymmetric matrices of odd size must have a real eigenvalue, and it must be 0. Consequently, antisymmetric matrices of odd size are singular and in particular, at $\gamma=0$, $\hat W$ is singular.
\end{remark}

Considering practical applications, where $\hat W$ is computed from data, we can expect $\hat W$ to include a small not anti-symmetric perturbation $P$, such as
\begin{equation}\label{eq:W_hat}
\hat W(\gamma) = W- W^\top + \gamma \textnormal{Id} + P,
\end{equation}
where the dependence on $\gamma$ is now made explicit and where we expect $P$ to be small in a suitable matrix norm. It is expected that a generic small perturbation is sufficient to make the degeneracy disappear. Then, pairs of eigenvalues would still cross the imaginary axes, but not at the same values of $\gamma$. This can generate a finite cascade of Hopf bifurcations for $\gamma$ close to 0. Additionally, a perturbation would ensure that the Jacobian is non-singular at all values of $\gamma$, including odd-dimensional problems at $\gamma =0$. In summary, we now have an ODE representing an RNN, such that one hyperparameter could be responsible for a cascade of Hopf bifurcations. With such initial understanding in mind, we delve deeper into the search for Hopf bifurcations, and describe them from a more general perspective.

It is worth remarking that the following two sections are presented in an abstract setting and can be applied to a variety of RNNs. The code presented in \cite{ODEinRNN} is mostly general, and the modifications to apply it to a new system only include the definition of the system and its derivatives.

We apply the rewriting \eqref{eq:Hopf} of the Hopf bifurcation to equation \eqref{eq:ODE_AsymRNN}, thus retrieving an algebraic Hopf problem of the form
\begin{align}\label{eq:zero_Hopf_RNN}
F(y) &= F(x, \gamma, v_r, v_i, \lambda_i) \\&= 
\begin{pmatrix}
\phi^\top v_r  \\ \phi^\top v_i  - 1\\
\tanh(\hat W(\gamma) x ) \\
\partial_x (\tanh(\hat W(\gamma) x )) v_r  + \lambda_i v_i \\
\partial_x (\tanh(\hat W(\gamma) x )) v_i  - \lambda_i v_r \\
\end{pmatrix}\nonumber\\& = 0 ,\nonumber
\end{align}
where $F : \mathbb{R}^{3n + 2} \rightarrow \mathbb{R}^{3n + 2}$. 
Following the discussion in Section \ref{sec:BfInRNN}, we focus on the zero of $\tanh(\hat W(\gamma) x )$  happening at $x = 0$ for all $\gamma$. This simplifies $$\partial_x (\tanh(\hat W(\gamma) x )) = \tanh'(\hat W(\gamma) x )\hat W(\gamma) = \hat W(\gamma).$$ The Jacobian then has $\lfloor n/2\rfloor$ pairs of complex eigenvalues, each pair crossing the imaginary axis at some value of $\gamma$, likely close to 0. We consider then each pair of eigenvalues $\lambda_j, j = 1,\dots, n/2$ as a function of $\gamma$. For each pair of eigenvalues, we search for the bifurcation value $\gamma_j$, such that $\lambda_j(\gamma_j)$ is purely imaginary. Then, we validate the existence of an algebraic Hopf bifurcation at $(x,\gamma)=(0, \gamma_j)$. We also study the other eigenvalues of $\hat W(\gamma_j)$ to ensure that no other eigenvalue is crossing the imaginary axis at the same value $\gamma_j$. After computing the first Lyapunov coefficient and ensuring that it is different from $0$, we have completed the validation of Hopf bifurcations.
Furthermore, it is possible to check the sign of the first Lyapunov coefficient to confirm the stability of the periodic orbit generated at the Hopf bifurcation. All run validations returned a positive Lyapunov exponent, determining that the bifurcating branch of periodic solutions is unstable.


\subsection{Hopf validation for RNNs}

We now provide the details for the concrete validation for our model problem. We define $\hat W(\gamma)$ as in \eqref{eq:W_hat}, where $W$ and $P$ are two random square matrices whose values are taken from a normal distribution. For our code, it is possible to give as input (to \texttt{asym\_RHS\_Hopf}) the seed for the random number generator and the amplitude of the perturbation, thus defining both $W$ and $P$.

For a wide array of problems, with dimensions ranging from 2 to 400, the code provided in \cite{ODEinRNN} numerically computes the existence of a solution to the algebraic Hopf problem \eqref{eq:Hopf}, validates the algebraic Hopf problem following the procedure in Appendix \ref{sec:RPA_details} and computes the validated first Lyapunov coefficient as in Appendix \ref{sec:FLC}. In Figure \ref{fig:transient_and_diverging}, some visual results are presented. Each plot represents (some coordinates of) an orbit close to the Hopf bifurcation.
In the first column of Figure \ref{fig:transient_and_diverging}, transient behaviour is presented. In all cases computed, all Lyapunov coefficients are positive, thus the periodic orbits created by the Hopf bifurcations are unstable. Still, it is usually possible to numerically shadow them for a short period of time. This is achieved by starting a forward integration at an approximate periodic orbit. The second column of Figure \ref{fig:transient_and_diverging} presents an orbit at the same parameter $\gamma$ close to the Hopf bifurcation, but the initial condition is chosen far away from the periodic orbit. This gives a graphical comparison between an orbit initially shadowing periodicity, and generic orbit at the same parameter value. The parameter value chosen is slightly larger than the largest $\gamma_j$.

The validation algorithm is very robust, and succeeds in validating most numerically found Hopf bifurcations. For example, in the 400 dimensional case presented, 200 Hopf bifurcations are found numerically and 198 are validated. Rarely, the random $\hat W$ created can be singular at a Hopf bifurcation. If this is the case, the validation of the algebraic Hopf problem fails. 

\begin{figure}
\centering
\begin{subfigure}[t]{.45\textwidth}
  \begin{overpic}[width=0.9\linewidth]{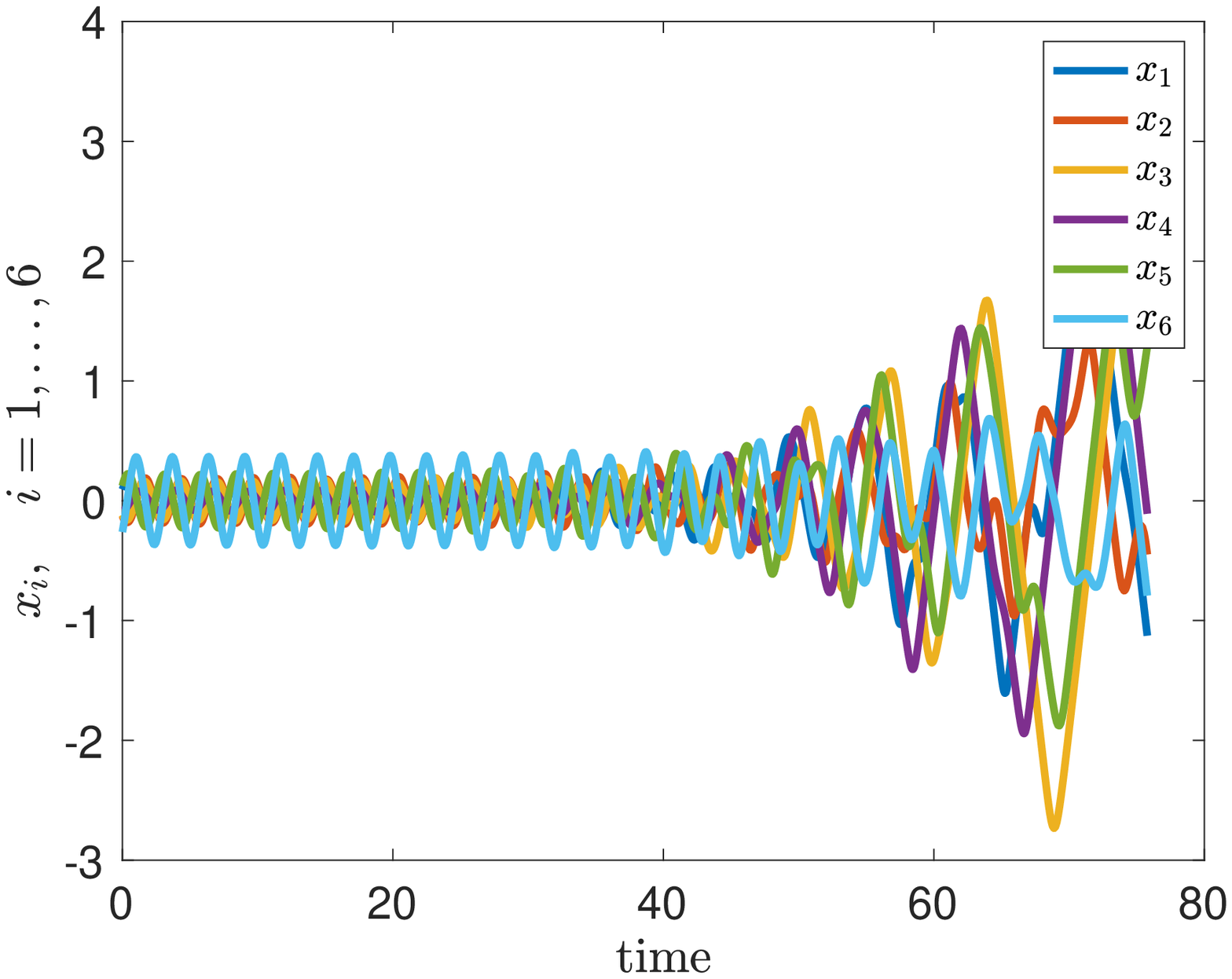}
\end{overpic}
  \caption{Transient dynamics in 6 dimensions, $\gamma = 0.025848$}
  \label{fig:dim6_transient}
\end{subfigure}%
\hspace{.4cm}
\begin{subfigure}[t]{.45\textwidth}
    \begin{overpic}[width=0.9\linewidth]{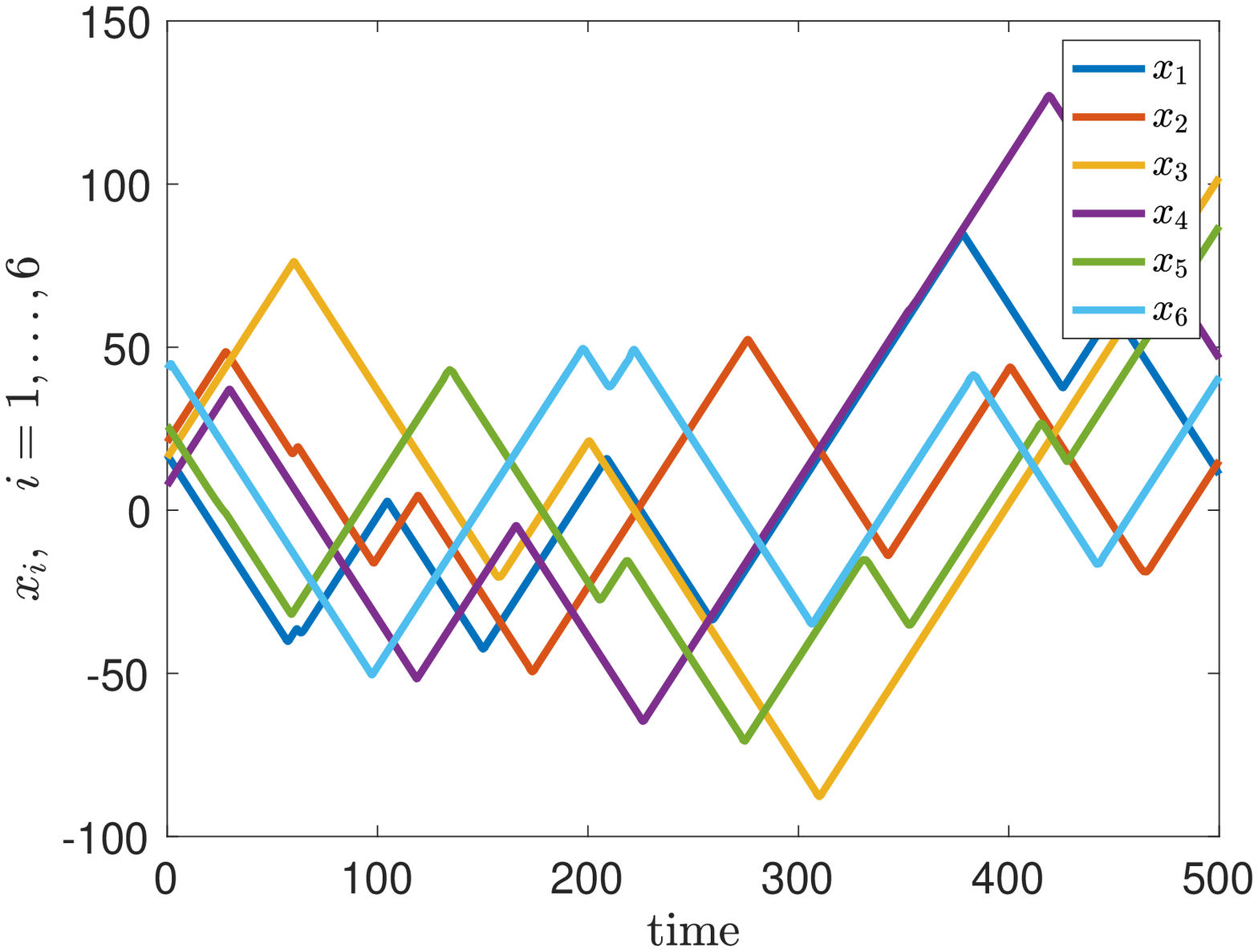}
\end{overpic}
  \caption{Another orbit in 6 dimensions for the same parameter value}
  \label{fig:dim6_divergent}
\end{subfigure}
\vspace{0.6cm}

\begin{subfigure}[t]{.45\textwidth}
  \begin{overpic}[width=0.9\linewidth]{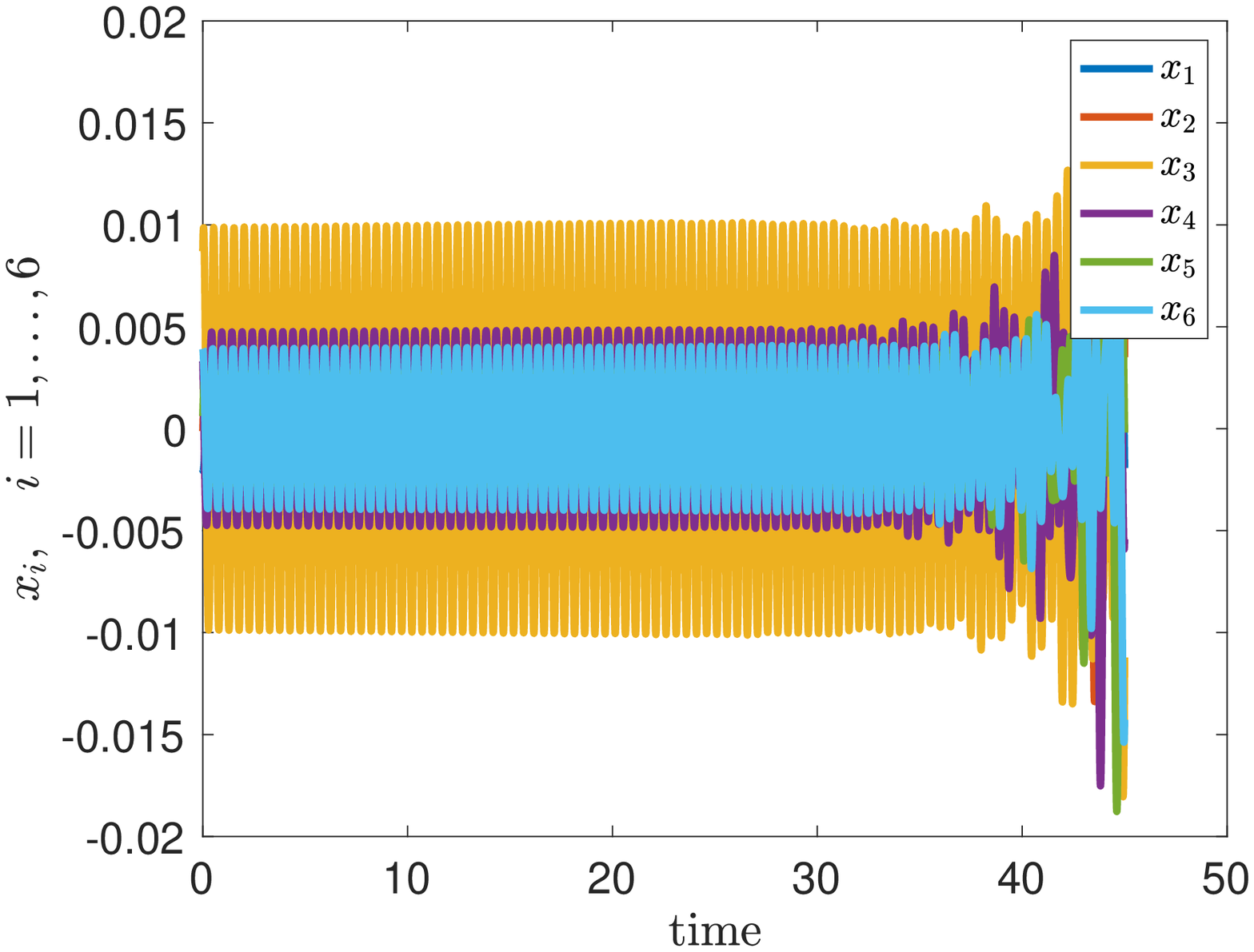}
\end{overpic}
  \caption{Transient dynamics in 50 dimensions, $\gamma = 0.193188$}
  \label{fig:dim50_transient}
\end{subfigure}%
\hspace{.4cm}
\begin{subfigure}[t]{.45\textwidth}
    \begin{overpic}[width=0.9\linewidth]{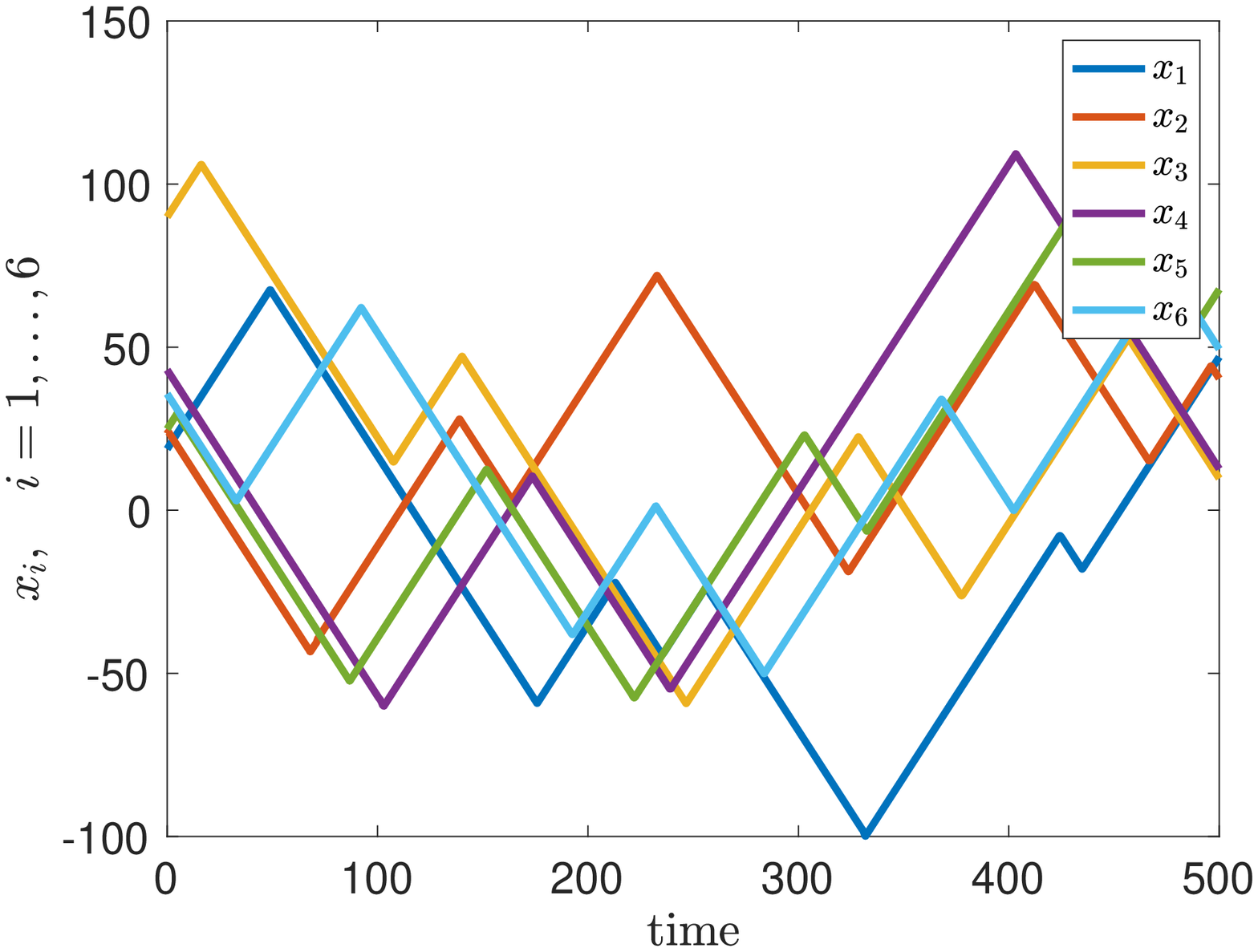}
\end{overpic}
  \caption{Another orbit in 50 dimensions for the same parameter value}
  \label{fig:dim50_divergent}
\end{subfigure}
\vspace{0.6cm}

\begin{subfigure}[t]{.45\textwidth}
  \begin{overpic}[width=0.9\linewidth]{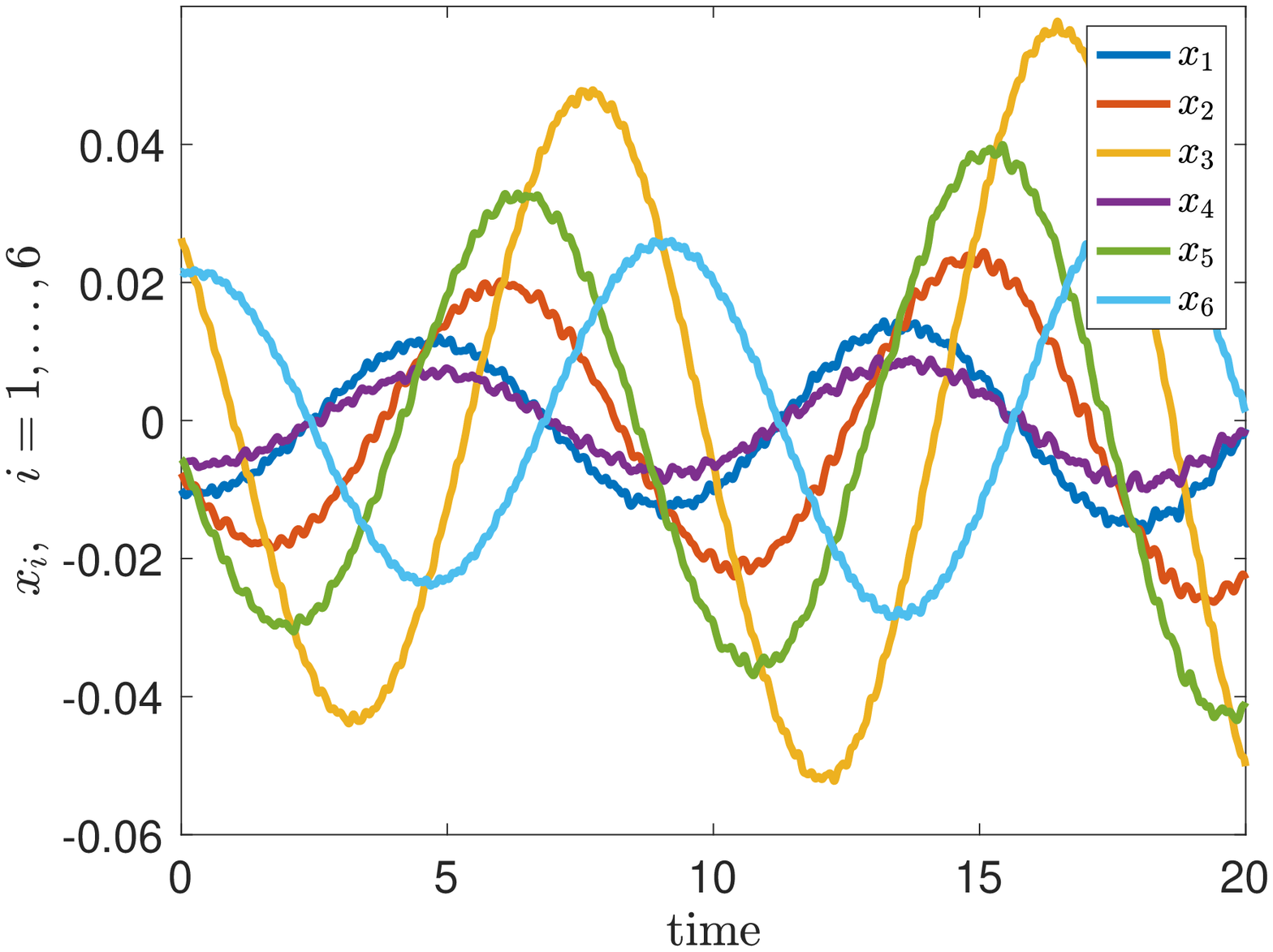}
\end{overpic}
  \caption{Transient dynamics in 400 dimensions, $\gamma = 0.201686$}
  \label{fig:dim400_transient}
\end{subfigure}%
\hspace{.4cm}
\begin{subfigure}[t]{.45\textwidth}
\begin{overpic}[width=0.9\linewidth]{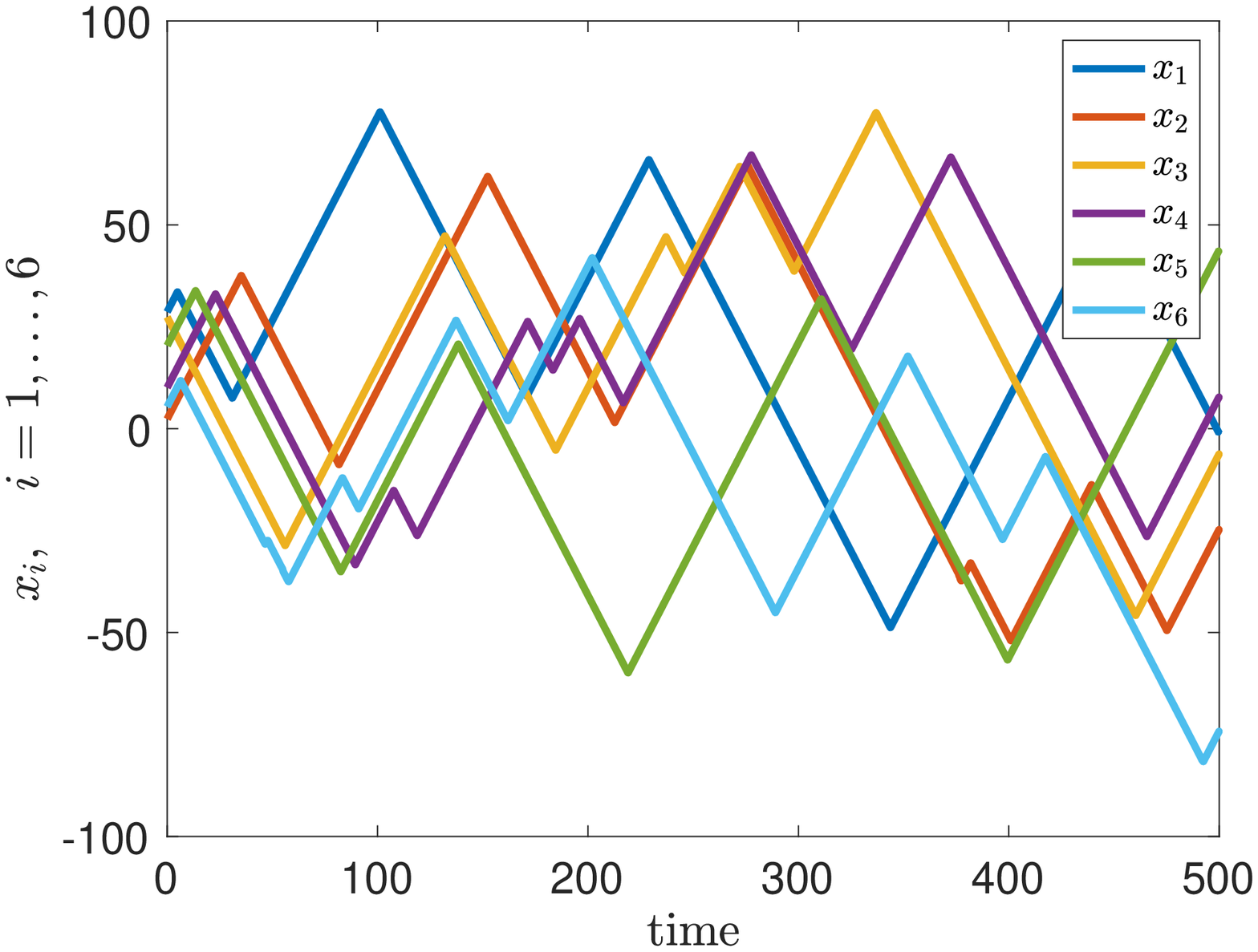}
\end{overpic}
  \caption{Another orbit in 400 dimensions for the same parameter value}
  \label{fig:dim400_divergent}
\end{subfigure}

\caption{For a fixed parameter $\gamma$ close to a Hopf bifurcation, transient dynamics is compared to a random orbit in 6, 50 and 400 dimensions. For these computations, the matrices $W$ and $P$ are set as randomly generated matrices, each element is i.i.d.~in $\mathcal{N}(0,1)$ and $\mathcal{N}(0,0.1)$ respectively. The seed for the random number generator is set to 80.
In the interest of clarity, only the first 6 coordinates are plotted with respect to time.
}
\label{fig:transient_and_diverging}
\end{figure}

\noindent
point before diverging.

All figures can be created running \texttt{figure\_generation.m}.



\subsection{Numerical comparison}

\begin{figure}
\centering
\begin{subfigure}[t]{.45\textwidth}
\includegraphics[width=0.9\linewidth]{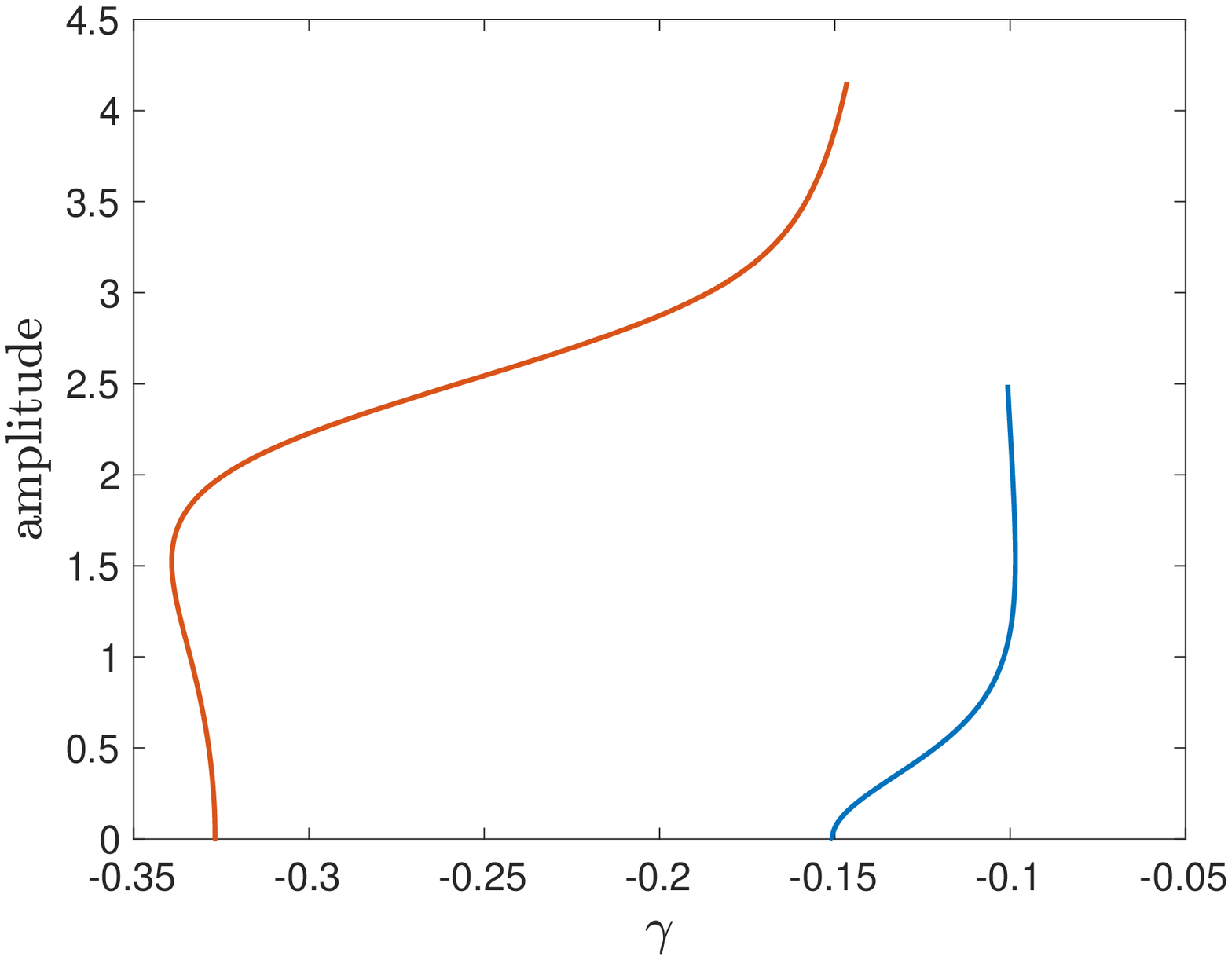}
\caption{Bifurcation diagram associated with ODE \eqref{eq:4D_bif_diag}.}
\label{fig:bif_diag4D}
\end{subfigure}%
\hspace{.4cm}
\begin{subfigure}[t]{.45\textwidth}
\includegraphics[width=0.9\linewidth]{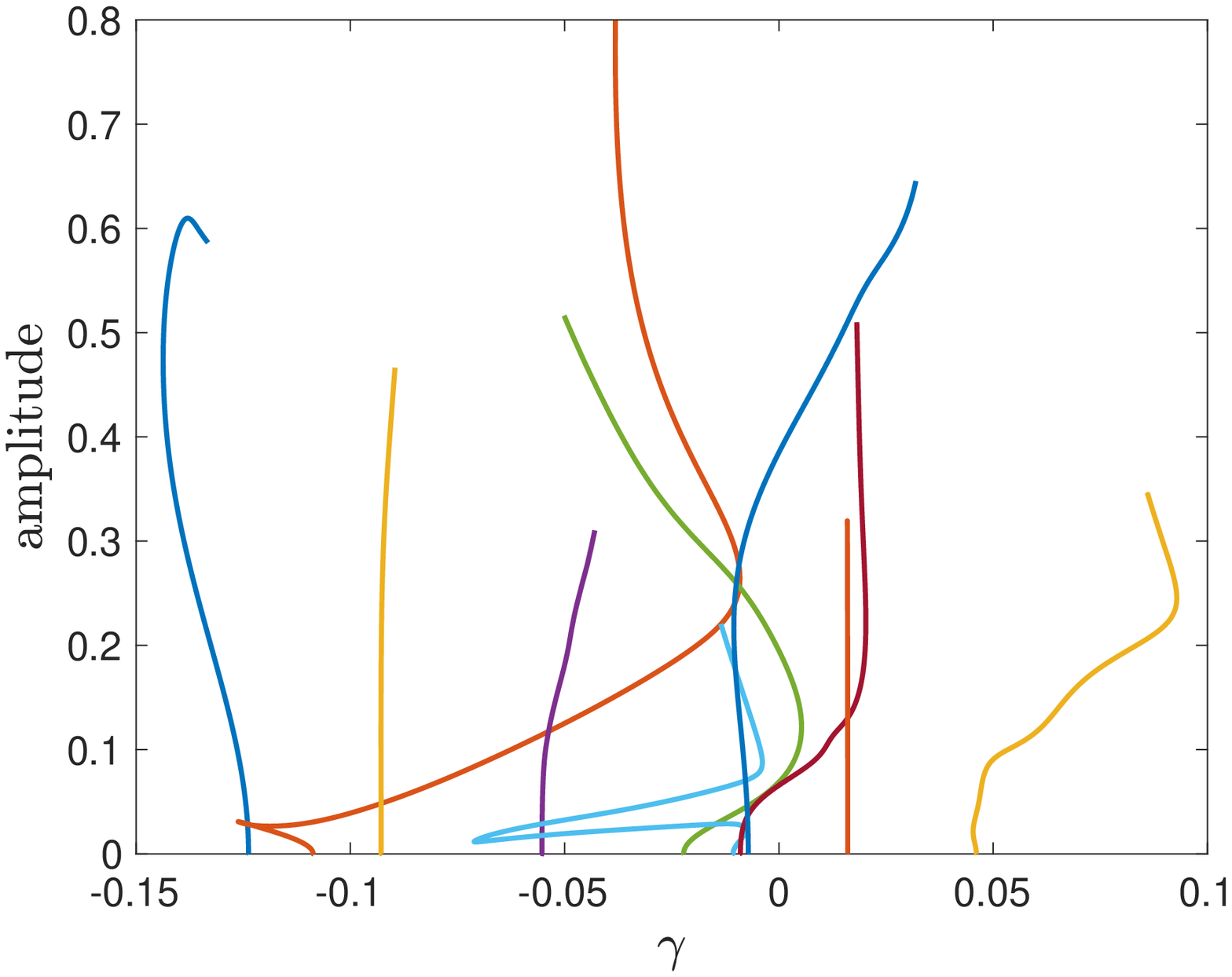}
\caption{An example of bifurcation diagram for a 20 dimensional system, with  randomly generated matrices $W$ and $P$ and with seed 80.}
\label{fig:bif_diag20D}
\end{subfigure}%

\caption{Bifurcation diagrams. In these figures, the Hopf bifurcations are validated, while the periodic orbits are continued numerically from the bifurcation using MatCont. Here, we plot the parameter $\gamma$ versus the amplitude of the periodic orbits.}

\end{figure}

To study more in depth the creation of new periodic orbits, we also wanted to understand, whether the branches of periodic orbits match purely numerical results (i.e.~without validation) that can be obtained with numerical continuation software. We have used MatCont~\cite{dhooge2003matcont}, which provides a nice cross-benchmark. For this comparison, we started with the system
\begin{equation}\label{eq:4D_bif_diag}
    \begin{cases}
x' = -\tanh((0.1\gamma+0.0929)x + 1.4109y \\\qquad -0.6359z -1.6482w) ,\\
y' = -\tanh(-1.3993x +(0.1\gamma+-0.0672)y \\\qquad +0.8243 z+0.7872w),\\
z' = -\tanh(0.7769x -0.7604y +(0.1\gamma\\\qquad + 0.0325) z+1.8087w) , \\
w' = -\tanh(1.4191x -0.8182y -1.7241z\\\qquad+(0.1\gamma+0.0373)w) ,
    \end{cases}
\end{equation}
where the initial matrices where created at random, and the Hopf bifurcations were validated before the rescaling of the parameter $\gamma$. The rescaling was implemented to increasingly separate the two Hopf bifurcations. 
For these parameters, the bifurcation diagram is depicted in Figure \ref{fig:bif_diag4D}. The Hopf bifurcations computed with MatCont agree with the rigorously validated ones. The global continuation of the periodic orbit branches is purely numerical, and can be achieved with MatCont~\cite{dhooge2003matcont}. Validation of periodic orbits starting from a Hopf bifurcation is a topic presented in \cite{queirolo2021Hopf}. The full presentation of these techniques is not within the scope of this article, and we refer the interested read to \cite{queirolo2021Hopf} for an overview of validation of periodic orbits generated from Hopf bifurcations in polynomial vector fields. In our case, the vector field is non-polynomial and automatic validation techniques need to be applied, such as the ones presented in \cite{lessard2016automatic} and \cite{groothedde2017parameterization}. 

In a similar way, larger ODE systems can be studied, where the bottleneck resides mainly in computational time. In this way, a 20 dimensional system has been studied, the code presented in {\texttt{validation\_2\_matcont.m}}. It was possible to validate 10 Hopf bifurcations, as expected from a 20 dimensional system, but MatCont was able to find only 9 of them, mainly due to bifurcations being very close to one-another. If properly initialised, MatCont could recognise the validated bifurcations as Hopf bifurcations. The full bifurcation diagram is presented in Figure \ref{fig:bif_diag20D}, it has been created using MatCont, and as such it is not validated globally but only locally near the bifurcation points. 

\section{General example}\label{sec:examples}

The assumption \eqref{eq:Wstructure} presented in the previous section allows us to analytically support the intuition for the existence of Hopf bifurcations in antisymmetric RNNs, but it is not a restriction neither in the analysis nor in the code here presented. In this section, we present a variety of RNNs structures and the associated proofs. The examples presented in this section drop the constraint of the antisymmetry of the weight matrix and consider different smooth activation functions, as well as different parameters for the Hopf bifurcation. Our final example deals with a multiple layer RNN, where then the right hand side of the Neural ODE is the subsequent application of the usual layer structure.

\subsection{Perturb a single diagonal element}

\begin{figure}
\centering
\includegraphics[width=0.5\linewidth]{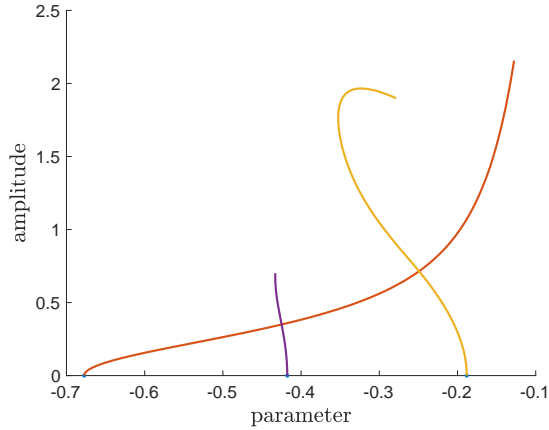}
\caption{Bifurcation diagram associated with Equation \eqref{eq:diag_par}, with randomly generated weight matrices $W$ and $P$ with seed 90.}
\label{fig:diag_par}
\end{figure}

In this example, we maintain the general structure presented in \eqref{eq:ODE_AsymRNN}, but we change the parameter, thus redefining
\begin{equation}\label{eq:diag_par}
\hat W (a) = W - W^T + \gamma Id + a D,
\end{equation}
where $D$ is a matrix with a unique non-zero element on the diagonal. For the purpose of this example, we chose a 6-dimensional system and we chose $D_{6,6} = 1$. We could then find and validate 3 Hopf bifurcations, creating the bifurcation diagram shown in Figure \ref{fig:diag_par}.

\subsection{Perturb an off-diagonal element}

\begin{figure}
\centering
\includegraphics[width=0.5\linewidth]{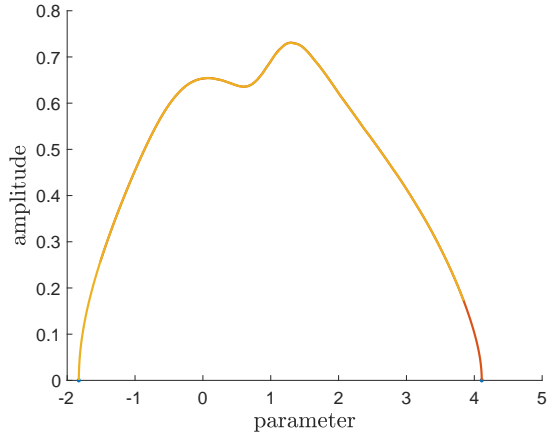}
\caption{Bifurcation diagram associated with Equation \eqref{eq:offdiag_par}. We found many systems sustaining periodic branches for larger intervals of parameter space in comparison to the antisymmetric case.}
\label{fig:offdiag_par}
\end{figure}

In a similar way, we can consider to modify a single off-diagonal element, thus breaking the antisymmetry of the system, by considering
\begin{equation}\label{eq:offdiag_par}
\hat W (a) = W - W^T + \gamma Id + a F,
\end{equation}
where $F$ is a matrix with a unique non-zero element not on the diagonal. For the purpose of this example, we chose a 6-dimensional system and we chose $F_{4,3} = 1$. It would be also possible to consider $F$ itself an anti-symmetric matrix, thus preserving the original structure. Such additional example is included in the code. We then find and validate one Hopf bifurcation, creating the bifurcation diagram shown in Figure \ref{fig:offdiag_par}.

\subsection{Multiple layers}

\begin{figure}
\centering
\includegraphics[width=0.5\linewidth]{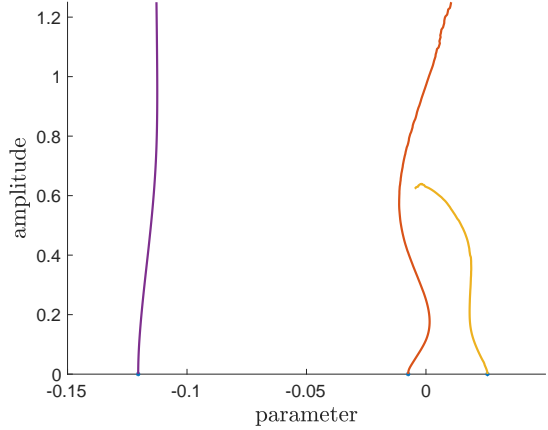}
\caption{Bifurcation diagram associated with Equation \eqref{eq:multilayer}.}
\label{fig:multilayer}
\end{figure}

As a final example, we consider a system with multiple layers. We define the new neural ODE as
\begin{equation}\label{eq:multilayer}
\begin{cases}
    \dot x = x_n,\\
    x_{i+1} = \sigma (W_i x_i), \quad i = 0,\dots n-1,\\
    x_0 = x,
\end{cases}
\end{equation}
where $\sigma$ is the activation function, in this cases chosen to be the hyperbolic tangent, and $W_i$ are matrices of weights depending on a parameter, each satisfying Equation \eqref{eq:W_hat}. We remark that these weights all depend on the same parameter $\gamma$, but this choice has been motivated by the simplicity of the exposition. For this example, we chose a 6-dimensional system with 3 layers. We remark how the complexity of the right hand side makes the computations for systems of this type scale by the product between dimensions and layers. As we can see in Figure \ref{fig:multilayer}, also in this situation we could find and validate Hopf bifurcations, and generate a (non-validated) bifurcation diagram. The behaviour we observe here is similar to the one we already encountered in the single layer system.

\section{Outlook}
\label{sec:outlook}

In this work, we have initiated the automated rigorous study of validation methods for nonlinear neural network dynamics. In a case study example, we have demonstrated the power of validation techniques to prove the existence of large numbers of periodic orbits, which can lead to complex oscillatory transients in the studied class of recurrent neural networks.\medskip

We briefly comment on potential, even more general, applications for the rigorous validation paradigm. First, we could aim to study the learning process as a dynamical system, or even couple information propagation on the network with learning dynamics of the network (i.e., of its weights). The major obstacle in this context is not a primarily a conceptual one but a computational obstacle as the dynamical system for the weights grows quadratically in the dimension of the nodes. Second, we aim to validate in future work even more complex structures beyond periodic orbits and also tackle more global bifurcation curves. In particular, the main idea to combine rigorous computer validation to mathematically prove results about neural network dynamics is broadly applicable.\medskip  

\textbf{Acknowledgements:} The reformulation \eqref{eq:Hopf} of algebraic Hopf bifurcation comes from a collaboration of Elena Queirolo with Jean-Philippe Lessard. CK \thanks{ckuehn@ma.tum.de} and EQ \thanks{elena.queirolo@tum.de}  would like to thank the VolkswagenStiftung for support via a Lichtenberg Professorship. EQ would like to thank the German Science Foundation (Deutsche Forschungsgemeinschaft, DFG) for support via a Walter-Benjamin Grant.

\appendix

\section{Proving algebraic Hopf bifurcations}
\label{sec:RPA_details}

Given a vector $\phi\in\mathbb{R}^n$, the algebraic Hopf problem \eqref{eq:Hopf} is a zero finding problem
\begin{align}\label{eq:zero_finding}
F: &X = \mathbb{R}\times\mathbb{R}\times\mathbb{R}^n\times\mathbb{R}^n\times\mathbb{R}^n \\&\quad \rightarrow Y = \mathbb{R}\times\mathbb{R}\times\mathbb{R}^n\times\mathbb{R}^n\times\mathbb{R}^n \nonumber\\
&(x, \gamma, \lambda_i, v_r, v_i) \nonumber\\&\quad\mapsto 
\begin{pmatrix}
&\phi^\top v_r  \\ 
&\phi^\top v_i  - 1\\
&f(x, \gamma) \\
&\textnormal{D}_x f(x, \gamma) v_r  +\lambda_i v_i \\
&\textnormal{D}_x f(x, \gamma) v_i - \lambda_i v_r\\
\end{pmatrix}\nonumber
\end{align}
where the space $X = Y$ is endowed with the Euclidean norm, thus making it a Banach space. In Section \ref{sec:BfInRNN}, we are interested in 
\begin{equation}\label{eq:ODE_RNN}
f(x, \gamma) = \tanh( \hat W(\gamma) x),
\end{equation}
but other definitions can be used with equal success. For ease of exposition we will be using this choice of $f$ in the following derivations.
For this case, the Jacobian is  defined in \eqref{eq:Jacobian}, and the algebraic Hopf zero-finding problem can be explicitly written as
\begin{align*}
&F(y ) \\&= \begin{pmatrix}&\phi^\top v_r  \\ 
&\phi^\top v_i  - 1\\
&\tanh( \hat W(\gamma) x) \\
&(\textnormal{Id} - \diag(\tanh(\hat W(\gamma) x)^2)) \hat W(\gamma) v_r  +\lambda_i v_i \\
&(\textnormal{Id} - \diag(\tanh(\hat W(\gamma) x)^2)) \hat W(\gamma) v_i - \lambda_i v_r\\\end{pmatrix}\\& = 0
\end{align*}
Let us assume we have a numerical solution $\hat y = (\hat x, \hat \gamma,\hat \lambda_i, \hat v_r, \hat v_i)$ to this zero finding problem. 

\begin{remark}\label{rem:initial}
Usually, iterative zero finding algorithms require a starting point $y_0$ in the neighborhood of a solution. For a general $f$, we can construct a rough approximation of the solution of \eqref{eq:zero_finding} by fixing $\gamma= \gamma_0$ and  finding $x_0$ numerically solving $f(x_0, \gamma_0) = 0$.  Then, let $(\lambda_j, v_j), j = 1,\dots,n$ be all the eigenpairs associated to $\textnormal{D}_xf(x_0,\gamma_0)$. A reasonable starting point is  $y_0=( x_0,  \gamma_0, \operatorname{Imag}(\lambda_0), v_0)$, where $\lambda_0$ is the eigenvalue closest to the imaginary axes and $v_0$ its associated eigenvector. Furthermore, the real part of $\lambda_0$ is a lower bound of the error. If such quantity is too large, the zero finding algorithm might diverge, and the initial choice of $\gamma_0$ should be reconsidered.  Once a first approximation $y_0$ is found, the Newton's method can be used to return a sharper approximation $\hat y$. For the specific case \eqref{eq:ODE_RNN}, using the analytical knowledge of the problem we can fix $\gamma_0 = 0$ and $x=0$. Then, any eigenpair can be chosen, since all of them cross the imaginary axes for an appropriate value of $\gamma$ close to $\gamma_0=0$.
\end{remark}

With a numerical solution in hand, the radii polynomial approach is applied. We define  $ A$ as a numerical inverse to $\textnormal{D}F(\hat y)$. Note how $A$ is computed and stored numerically and does not need to be validated. Let $T(y) \bydef y - AF(y)$, as in \eqref{eq:contraction_mapping}. We then compute
$$
 \|T(\hat y)\|_X = \| AF(\hat y)\|_X \leq Y.
$$
By using Intlab on all computations, and all the computations being finite, the $Y$  bound can be computed directly. The $Z$ bound is more complicated and first requires a splitting. We use
\begin{align*}
&\sup_{b,c\in B_1(0)\subset X}\| \textnormal{D}T(\hat y +rb)rc \| _X
\\&\quad= \sup_{b\in B_1(0)\subset X} \| \textnormal{Id} - A \textnormal{D}F(\hat y + rb ) \|_{B(X,X)} r\\&\quad
\leq \| \textnormal{Id} - A \textnormal{D}F(\hat y) \|r_{B(X,X)} +\\&\quad \sup_{b,c\in B_1(0)\subset X}\|A\textnormal{D}F(\hat y) - A\textnormal{D}F(\hat y + rb) \|_{B(X,X)}r
\end{align*}
We then define
$$
Z_1(r) = Z_1 r \geq \| \textnormal{Id} - A \textnormal{D}F(\hat y) \|r_{B(X,X)},
$$
that can again be computed directly. For the second term, we apply the mean value theorem and get
\begin{align*}
\sup_{b\in B_1(0)\subset X}\|A\textnormal{D}F(\hat y) - A\textnormal{D}F(\hat y + rb) \|_{B(X,X)}r \\ \leq \sup_{b,c\in B_1(0)\subset X}\|AD^2F(\hat y + rb)rc \|_{B(X,X)}r.
\end{align*}
A bound of this type can be achieved by computing the second derivative on the \emph{interval} $\hat y \pm R$, where $R$ in an \emph{a priori} upper bound of the validation radius $r^*$. The interval notation is considered element-wise in all components of $\hat y$. In the same way, the vector $c$ is replaced by the vector $\textbf{1}\bydef 0\pm 1$, the vector having the interval $[-1,+1]$ in all coordinates. These modifications give
$$
Z_2(r, R) = Z_2r^2  \geq \| A\textnormal{D}^2F(\hat y \pm R) \textbf{1} \|_{B(X,X)}r^2,
$$
computed using interval arithmetic. Having constructed all the bounds, we need to find a value of $r$ such that the radii polynomial $p(r) = Y + (Z_1r + Z_2r^2) - r$ is negative. This is an explicit computation since the radii polynomial is second order in $r$, giving
$$
r_\pm^* = \frac{ 1 - Z_1 \pm \sqrt{(Z_1-1)^2 - 4Y Z_2}
}{
2 Z_2
}.
$$
If such $r^*_\pm$ exist, then the validation succeeded, and, for any $r$ in the interval  $[r^*_-,r^*_+]$, the ball centered at $\hat y$ of radius $r$ contains a unique solution. In our case, this finishes the proof that an algebraic Hopf bifurcation is taking place at most $r^*_-$ away from $(\hat x, \hat \gamma)$.

\section{Proving non-degeneracy of Hopf bifurcations} 
\label{sec:FLC}

This section is based upon standard results on Lyapunov coefficients presented in~\cite[Chapter~3]{kuznetsov2013elements} but we feel it is useful to explain in more detail, how non-degeneracy of Hopf is encoded within numerical validation techniques. A Hopf bifurcation $(x_\star, \gamma_\star)$ of the ODE $x' = f(x, \gamma)$ is non-degenerate if there is a unique pair of imaginary eigenvalues of the Jacobian $\textnormal{D}_xf(x_\star, \gamma_\star)$ and the first Lyapunov coefficient is non-zero. For the proof of this statement we refer to  \cite{kuznetsov2013elements}. In this section, we define the first Lyapunov coefficient constructively, such that its validated computation can conclude the proof of the non-degeneracy of a Hopf bifurcation. Let $x'=f(x, \gamma)$ be a parameter-dependent ODE, as in Section \ref{sec:finding_Hopf}, with $x\in\mathbb{R}^n$, $\gamma\in\mathbb{R}$, and let $(x_\star, \gamma_\star)$ be an algebraic Hopf bifurcation. Let $J$ be the Jacobian $\textnormal{D}_xf(x_\star, \gamma_\star)$, having a unique pair of purely imaginary eigenvalues. Let $\imag\lambda$ be the positive imaginary eigenvalue of $J$. Then, let $v$ be the eigenvector of $J$ associated to $\imag\lambda$, and let $w$ be the eigenvector of $J^\top$ associated to $-\imag\lambda$. Notice how both $v$ and $w$ are defined up to a complex constant. We introduce the complex inner product 
$$
\langle x,y\rangle \bydef \bar{x}^\top y,
$$
and we request $\langle v, v\rangle = 1$ and $\langle v,w\rangle = 1$. This rescaling is not necessary to determine the sign of the first Lyapunov coefficient, and could be skipped, as long as $\langle v,w\rangle >0$. Then, the first Lyapunov coefficient is defined as
\begin{align*}
l_1 = \frac{1}{2\lambda^2} \operatorname{Real}(
\imag
\langle w,\textnormal{D}^2f(x_\star,\gamma_\star)vv\rangle\\ \langle w,\textnormal{D}^2f(x_\star,\gamma_\star)v\bar v\rangle+\\
\lambda \langle w,\textnormal{D}^3f(x_\star,\gamma_\star)vv\bar v\rangle
).
\end{align*}
Each derivative of order $k$ is considered as an operator acting on $k$ elements of $\mathbb{R}^n$. Having built in the previous Section \ref{sec:RPA_details} an error bound of $x_\star$ and $\gamma_\star $, we use the Intlab eigenpair validation functionality to build a validated interval for $\lambda$, $v$ and $w$ and then we compute $l_1$ explicitly. If the validated interval of existence for $l_1$ does not contain 0, the proof of the non-degeneracy of the Hopf bifurcation is completed. In the examples provided, the majority of non-degeneracy validations were successful.


\begin{thebibliography}{10}

\bibitem{Alpaydin}
E.~Alpaydin.
\newblock {\em Introduction to Machine Learning}.
\newblock MIT Press, 2020.

\bibitem{AmbrosioGigliSavare}
L.~Ambrosio, N.~Gigli, and G.~Savar{\'e}.
\newblock {\em Gradient Flows: In Metric Spaces and in the Space of Probability
  Measures}.
\newblock {Birkh\"auser}, 2006.

\bibitem{ArioliKoch1}
G.~Arioli and H.~Koch.
\newblock Computer-assisted methods for the study of stationary solutions in
  dissipative systems, applied to the {Kuramoto–Sivashinski} equation.
\newblock {\em Arch. Rat. Mech. Anal.}, 197(3):1033--1051, 2010.

\bibitem{breden2013global}
Maxime Breden, Jean-Philippe Lessard, and Matthieu Vanicat.
\newblock Global bifurcation diagrams of steady states of systems of pdes via
  rigorous numerics: a 3-component reaction-diffusion system.
\newblock {\em Acta applicandae mathematicae}, 128(1):113--152, 2013.

\bibitem{chang2019antisymmetricrnn}
Bo~Chang, Minmin Chen, Eldad Haber, and Ed~H Chi.
\newblock Antisymmetricrnn: A dynamical system view on recurrent neural
  networks.
\newblock {\em arXiv preprint arXiv:1902.09689}, 2019.

\bibitem{chen2018neural}
R.T.Q. Chen, Yulia Rubanova, Jesse Bettencourt, and David Duvenaud.
\newblock Neural ordinary differential equations.
\newblock {\em arXiv preprint arXiv:1806.07366}, 2018.

\bibitem{dhooge2003matcont}
Annick Dhooge, Willy Govaerts, and Yu~A Kuznetsov.
\newblock Matcont: a matlab package for numerical bifurcation analysis of odes.
\newblock {\em ACM Transactions on Mathematical Software (TOMS)},
  29(2):141--164, 2003.

\bibitem{E5}
W.~E.
\newblock A proposal on machine learning via dynamical systems.
\newblock {\em Commun. Math. Stat.}, 5(1):1--11, 2017.

\bibitem{Silveretal}
D.~Silver et~al.
\newblock Mastering the game of {Go} with deep neural networks and tree search.
\newblock {\em Nature}, 529(7587):484--489, 2016.

\bibitem{FarmerPackardPerelson}
J.D. Farmer, N.H. Packard, and A.S. Perelson.
\newblock The immune system, adaptation, and machine learning.
\newblock {\em Phys. D}, 22(1):187--204, 1986.

\bibitem{GaltierWainrib}
M.N. Galtier and G.~Wainrib.
\newblock Multiscale analysis of slow-fast neuronal learning models with noise.
\newblock {\em J. Math. Neurosci.}, 2:13, 2012.

\bibitem{groothedde2017parameterization}
Chris~M Groothedde and JD~Mireles James.
\newblock Parameterization method for unstable manifolds of delay differential
  equations.
\newblock {\em Journal of Computational Dynamics}, 4(1\&2):21, 2017.

\bibitem{GH}
J.~Guckenheimer and P.~Holmes.
\newblock {\em Nonlinear Oscillations, Dynamical Systems, and Bifurcations of
  Vector Fields}.
\newblock Springer, New York, NY, 1983.

\bibitem{hungria2016rigorous}
Allan Hungria, Jean-Philippe Lessard, and Jason~D Mireles~James.
\newblock Rigorous numerics for analytic solutions of differential equations:
  the radii polynomial approach.
\newblock {\em Mathematics of Computation}, 85(299):1427--1459, 2016.

\bibitem{Ilyashenko2}
Yu. Ilyashenko.
\newblock Centennial history of {Hilbert's 16th} problem.
\newblock {\em Bull. Amer. Math. Soc.}, 39(3):301--354, 2002.

\bibitem{JordanMitchell}
M.I. Jordan and T.M. Mitchell.
\newblock Machine learning: Trends, perspectives, and prospects.
\newblock {\em Science}, 349(6245):255--260, 2015.

\bibitem{Kapelaetal}
T.~Kapela, M.~Mrozek, D.~Wilczak, and P.~Zgliczy{\'n}ski.
\newblock {CAPD:: DynSys:} a flexible {C++} toolbox for rigorous numerical
  analysis of dynamical systems.
\newblock {\em Commun. Nonl. Sci. Numer. Simul.}, 101:105578, 2021.

\bibitem{KuehnBook}
C.~Kuehn.
\newblock {\em Multiple Time Scale Dynamics}.
\newblock Springer, 2015.

\bibitem{ODEinRNN}
Christian Kuehn and Elena Queirolo.
\newblock Code for "computer validation of neural network dynamics: A first
  case study".
\newblock \url{https://github.com/elenaquei/RNNs/releases/tag/V1}, 2023.

\bibitem{kuznetsov2013elements}
Yuri~A Kuznetsov.
\newblock {\em Elements of applied bifurcation theory}, volume 112.
\newblock Springer Science \& Business Media, 2013.

\bibitem{LeCunBengioHinton}
Y.~LeCun, Y.~Bengio, and G.~Hinton.
\newblock Deep learning.
\newblock {\em Nature}, 521(7553):436--444, 2015.

\bibitem{lessard2016automatic}
Jean-Philippe Lessard, JD~Mireles James, and Julian Ransford.
\newblock Automatic differentiation for fourier series and the radii polynomial
  approach.
\newblock {\em Physica D: Nonlinear Phenomena}, 334:174--186, 2016.

\bibitem{MeiMontanariNguyen}
S.~Mei, A.~Montanari, and P.A. Nguyen.
\newblock A mean field view of the landscape of two-layer neural networks.
\newblock {\em Proc. Natl. Acad. Sci. USA}, 115(33):E7665--E7671, 2018.

\bibitem{DeMeloVanStrien}
W.~De Melo and S.~Van Strien.
\newblock {\em One-dimensional Dynamics}.
\newblock Springer, 2012.

\bibitem{Murphy}
K.P. Murphy.
\newblock {\em Machine Learning: a probabilistic perspective}.
\newblock MIT press, 2012.

\bibitem{Pinkus}
A.~Pinkus.
\newblock Approximation theory of the mlp model in neural networks.
\newblock {\em Acta Numer.}, 8:143--195, 1999.

\bibitem{rump1999intlab}
Siegfried~M Rump.
\newblock Intlab—interval laboratory.
\newblock In {\em Developments in reliable computing}, pages 77--104. Springer,
  1999.

\bibitem{DayLessardMischaikow}
J.P.~Lessard S.~Day and K.~Mischaikow.
\newblock Validated continuation for equilibria of {PDEs}.
\newblock {\em SIAM J. Numer. Anal.}, 45(4):1398--1424, 2007.

\bibitem{Sanger}
T.D. Sanger.
\newblock Optimal unsupervised learning in a single-layer linear feedforward
  neural network.
\newblock {\em Neural Networks}, 2(6):459--473, 1989.

\bibitem{ScarselloTsoi}
F.~Scarselli and A.C. Tsoi.
\newblock Universal approximation using feedforward neural networks: A survey
  of some existing methods, and some new results.
\newblock {\em Neural Networks}, 11(1):15--37, 1998.

\bibitem{Schmidhuber}
J.~Schmidhuber.
\newblock Deep learning in neural networks: An overview.
\newblock {\em Neural Networks}, 61:85--117, 2015.

\bibitem{SchoelkopfSmola}
B.~Sch{\"o}lkopf and A.J. Smola.
\newblock {\em Learning with Kernels: support vector machines, regularization,
  optimization, and beyond}.
\newblock MIT Press, 2001.

\bibitem{Tanaka}
T.~Tanaka.
\newblock Mean-field theory of {Boltzmann} machine learning.
\newblock {\em Phys. Rev. E}, 58(2):2302, 1998.

\bibitem{Tucker}
W.~Tucker.
\newblock The {Lorenz} attractor exists.
\newblock {\em C.R. Acad. Sci. Paris}, 328:1197--1202, 1999.

\bibitem{queirolo2021Hopf}
Jan~Bouwe Van~den Berg, Jean-Philippe Lessard, and Elena Queirolo.
\newblock Rigorous verification of hopf bifurcations via desingularization and
  continuation.
\newblock {\em SIAM Journal on Applied Dynamical Systems}, 20(2):573--607,
  2021.

\bibitem{vandenBergJamesReinhardt}
J.B. van~den Berg, J.D.M. James, and C.~Reinhardt.
\newblock Computing (un)stable manifolds with validated error bounds:
  non-resonant and resonant spectra.
\newblock {\em J. Nonlin. Sci.}, 26(4):1055--1095, 2016.

\bibitem{WangTengPerdikaris}
S.~Wang, V.~Teng, and P.~Perdikaris.
\newblock Understanding and mitigating gradient flow pathologies in
  physics-informed neural networks.
\newblock {\em SIAM J. Sci. Comput.}, 43(5):A3055--A3081, 2021.

\bibitem{yan2019robustness}
Hanshu Yan, Jiawei Du, Vincent~YF Tan, and Jiashi Feng.
\newblock On robustness of neural ordinary differential equations.
\newblock {\em arXiv preprint arXiv:1910.05513}, 2019.

\bibitem{zhang2019approximation}
Han Zhang, Xi~Gao, Jacob Unterman, and Tom Arodz.
\newblock Approximation capabilities of neural ordinary differential equations.
\newblock {\em arXiv preprint arXiv:1907.12998}, 2(4):3--1, 2019.

\end{thebibliography}
\end{document}